\documentclass[11pt]{amsart}
\usepackage{eucal}
\usepackage{amssymb}
\usepackage{epsfig}
\usepackage{amscd}
\usepackage{color}
\usepackage{enumerate}

\newcommand{\ra}{{\rightarrow}}

\renewcommand{\H}{\mathbb{H}}

\newcommand{\C}{\mathbb{C}}
\newcommand{\Z}{\mathbb{Z}}
\newcommand{\N}{\mathbb{N}}

\newcommand{\ie}{i.e.\ }
\newcommand{\Int}{\mathrm{Int}}

\def\pslc {\mathrm{PSL}(2,\C)}

\theoremstyle{plain}
\newtheorem{thm}{Theorem}[section]
\newtheorem{cor}[thm]{Corollary}

\newtheorem{lem}[thm]{Lemma}

\newtheorem{claim}[thm]{Claim}

\numberwithin{equation}{section}

\begin{document}
\title[Primitive stable representations]
        {Primitive stable representations of free Kleinian groups}
        \author[W. Jeon, I. Kim, and K. Ohshika with C. Lecuire]{Woojin Jeon, Inkang Kim, and Ken'ichi Ohshika\\
        \lowercase{with the collaboration of} Cyril Lecuire}
        \date{}
        \maketitle

\begin{abstract}
In this paper, we give a complete criterion
for a discrete faithful representation $\rho:F_n \ra \pslc$ to be primitive stable.
This will answer Minsky's conjectures about geometric conditions on $\H^3/\rho(F_n)$ regarding the primitive 
stability of $\rho$.
\end{abstract}
 \footnotetext[3]{The second
 author gratefully acknowledges the partial support of  NRF grant
(R01-2008-000-10052-0).}
\footnotetext[4]{The fourth author gratefully acknowledges the support of  Indo-French Research Grant 4301-1.}

\tableofcontents
\section{Introduction}
Let $F$ be a non-abelian free group of rank $n$.
For any group $G$, the automorphism group Aut$(F)$ acts on Hom$(F,G)=G^n$ by precomposition. This action projects down to the action of the outer automorphism group Out$(F)$ on the character variety ${\mathcal X}(F,G)$ which is defined as the geometric quotient of Hom$(F,G)$ by inner automorphisms of $G$.
When $G$ is $\pslc$, Minsky studied a dynamical decomposition of ${\mathcal X}(F,G)$ with respect to the Out$(F)$-action. 
Here a dynamical decomposition means decomposing ${\mathcal X}(F,G)$ in terms of proper discontinuity and ergodicity of the action. See \cite{Lu} for more information about this decomposition.  Minsky defined the set  $\mathcal{PS}(F)$ of {\it primitive stable characters}, and his main results are as follows (\cite{Min2}).
\begin{enumerate}
\item $\mathcal{PS}(F)$ is an open subset of ${\mathcal X}(F,\pslc)$, and Out$(F)$ acts on $\mathcal{PS}(F)$ properly discontinuously.
\item $\mathcal{PS}(F)$ is strictly larger than the set of Schottky characters.
\item For every proper free factor $A$ of $F$ and a primitive stable representation $\rho$, the restriction $\rho\vert_A$ is Schottky.
\end{enumerate}
Since the set of Schottky characters is known to be the interior of the set of discrete faithful characters by Sullivan \cite{Sul}, it follows that the dynamical decomposition of ${\mathcal X}(F,G)$ is different from the well-known geometric decomposition of ${\mathcal X}(F,G)$ in terms of being discrete faithful and having dense image.

When a primitive stable representation $\rho$ is discrete and faithful, $\H^3/\rho(F)$ becomes a hyperbolic 3-manifold which we call a {\it hyperbolic handlebody} and our main interest is finding a geometric condition on $\H^3/\rho(F)$ under which $\rho$ becomes primitive stable. 
In the same paper \cite{Min2}, Minsky conjectured that 
\begin{enumerate}
\item Every discrete faithful representation of $F$ without
parabolics is primitive stable.
\item A discrete faithful representation of $F$ is primitive stable
if and only if every component of its ending lamination is blocking.
\end{enumerate}
We shall see that a {\it disc-busting} minimal lamination is {\it blocking}.
For definitions, see  section \ref{pre}. In this paper, we shall prove the first conjecture 
and also give an answer to the second one.

\begin{thm} 
\label{main without parabolics}
If $\rho$ is a discrete faithful representation of $F$ without
parabolics then $\rho$ is primitive stable. 
\end{thm}
In this case, the ending lamination of $\H^3/\rho(F)$ is
necessarily connected and in the Masur domain.
\begin{thm}
\label{main with parabolics}
Let $\rho$ be a discrete, faithful and geometrically infinite
representation with parabolics such that the non-cuspidal part $M_0$ of $M=\H^3/\rho(F)$ is the union of the relative compact core $H$ and finitely many end neighbourhoods $E_i$ facing $S_i\subset \partial H$. 
Then the representation $\rho$ is primitive stable if and only if every parabolic curve is disc-busting, and every geometrically infinite end $E_i$ has the ending lamination $\lambda_i$ which is disc-busting on $\partial H$.
\end{thm}
Theorem \ref{main without parabolics} and the sufficiency part of Theorem \ref{main with parabolics} were announced in \cite{JK}. 
The idea of proof is following Minsky's construction of  a primitive stable representation which is not Schottky, and the main new ingredient comes from the recent result of Mj about Cannon-Thurston maps of free Kleinian groups (see \cite{Mj8}). 

\section{Preliminaries}\label{pre}
\subsection{Primitive stability} Let us recall that $F$ is a non-abelian free group of rank $n$. Let $\vee S^1$ be a bouquet of $n$ oriented circles, whose fundamental group is $F$ with a fixed generating
set $X=\{x_1,\cdots,x_n\}$. 
Then its universal cover $\widetilde{\vee S^1}$ can be identified with the Cayley graph $\Gamma_F$ of $F$ with respect to $X$, which is a tree with the canonical word metric. Following \cite{BFH}, the space of oriented lines $\widetilde {\mathcal B}(F)$ in $\Gamma_F$ can be identified with $(\partial_\infty F\times \partial_\infty F) \backslash \triangle$ where $\partial_\infty F$ denotes the Gromov boundary of the tree $\Gamma_F$ and $\triangle$ denotes the diagonal. 
The free group $F$ acts diagonally on $\widetilde {\mathcal B}(F)$ as the covering transformations. 
We denote the quotient space of $\widetilde {\mathcal B}(F)$ under this action by ${\mathcal B}(F)$, and call each element of ${\mathcal B}(F)$ also a line. 

For $w\in F$, if we let $\overline{w}$ be the bi-infinite periodic word determined by concatenating infinitely many copies of $w$, then it defines an $F$-invariant family of lines in
$\widetilde {\mathcal B}(F)$ and these lines are projected to a line in ${\mathcal B}(F)$ which can be identified with $\overline{w}$ modulo shift or the conjugacy class $[w]$.
An element of $F$ is called {\it primitive} if it can be a member of a free generating set and we let ${\mathcal P}(F)$ denote the subset of ${\mathcal B}(F)$ consisting of $\overline{w}$ for conjugacy classes $[w]$ of primitive elements, which is Out$(F)$-invariant.

Given a representation  $\rho : F \rightarrow\pslc$ and a base point $o\in
\H^3$, there is a unique $\rho$-equivariant map
$\tau_{\rho,o}:\Gamma_F \ra \H^3$ sending the origin $e$ of $\Gamma_F$
to $o$ and taking each edge to a geodesic segment \cite{Fl}. 
A representation $\rho : F \rightarrow \pslc$ is {\it primitive stable} if
there are constants $K,\delta$ such that
$\tau_{\rho,o}$ takes all lines (in $\widetilde {\mathcal B}(F)$) corresponding to $\mathcal{P}(F)$ to $(K,\delta)$-quasi-geodesics in $\H^3$.
This definition is independent of the choice of the base point $o\in \H^3$, which we can easily see by changing $\delta$.  
Since the primitive stability is also invariant under conjugacy, for simplifying arguments for checking primitive stability, we shall define a unique element of $\widetilde {\mathcal B}(F)$ corresponding to a cyclically reduced $w\in F$ as follows. Since $\Gamma_F$ is a tree,
there exists a unique oriented line $\widetilde w$ on $\Gamma_F$ passing through all the $w^k(e)$ for $k \in \Z$, where we regard $w$ as a covering transformation. 
Then clearly, the broken geodesic image $\tau_{\rho,o}(\widetilde w)$ passes through the base point $o$.

\subsection{Whitehead lemma} 
\label{sec: Whitehead}
We refer to \cite{CB, Th} for the basic theory of geodesic and measured laminations. 
The space of measured laminations on a hyperbolic surface is a completion of weighted simple geodesics, to which the geometric intersection number continuously extends. 
Recall that the {\it Masur domain}  \cite{Ma} of a handlebody $H$ consists of projective classes of measured laminations which have positive intersection number with every non-empty limit of weighted meridians of $H$.
A measured lamination $\lambda$ (or a simple closed curve)  is said to be {\it disc-busting} if there exists
$\eta>0$ such that for any essential disc  $A$,\ $i(\partial A,
\lambda)>\eta$. 
Otherwise $\lambda$ is called {\em disc-dodging}.


Let $\Delta=\{\delta_1,\cdots,\delta_n\}$ be a system of
compressing discs on a handlebody $H$ along which one can cut $H$
into a $3$-ball. 
We call such a system a {\em cut system}.
A free generating set $X$ of $\pi_1(H)=F$ is dual
to such a system.
The {\em Whitehead graph} $Wh(\lambda,\Delta)$ (of $\lambda$ with respect to $\Delta$) is defined as follows.
First, we isotope $\Delta$ so that its boundary consists of closed geodesics with respect to a fixed hyperbolic metric on $\partial H$(This removes all inessential intersections between $\lambda$ and $\partial \Delta$).
 By cutting $\partial H$ along $\Delta$, we get a planar surface with $2n$ boundary
components.
We label the two discs coming from cutting $H$ along $\delta_i$ as $\delta_i^+, \delta_i^-$ so that the oriented loop corresponding to a generator of $\pi_1(H)$ passes from $\delta_i^-$ to $\delta_i^+$. Then $\partial\delta_i^+$ and $\partial\delta_i^-$ are identified on $\partial H$. 
The vertices of the graph $Wh(\lambda,\Delta)$  are the boundary components of $\Delta$, and
two vertices are connected by an edge if and only if the corresponding boundary components are joined by an arc of
$\lambda\setminus \Delta$. 
Note that we can regard  the edges of $Wh(\lambda,\Delta)$ as embedded in $\partial H \setminus \Delta$, and we can replace the vertices with small circles $\partial\delta_i^{\pm}$.

Given a meridian or a system of meridians $m$ on $\partial H$, an {\it $m$-wave} is an arc $\alpha\subset\partial H$ such that $\alpha\cap m=\partial\alpha$ and that there is an arc $\beta\subset m$ such that $\alpha$ is homotopic to $\beta$ in $H$ but not on $\partial H$. A geodesic lamination $\lambda$ is said to be in {\it tight position} with respect to a cut system $\Delta$ (or in more general, to a disjoint union of compressing discs) if any $\partial\Delta$-wave intersects $\lambda$. Otal showed the following in Proposition 3.10 of  his thesis \cite{Ot}.

\begin{lem}
\label{tight}
Suppose that a geodesic lamination $\lambda$ is in a tight position with respect to a cut system $\Delta$.
Then $Wh(\lambda, \Delta)$ is connected and has no cut points.
\end{lem}

For disc-busting measured laminations, we can show the following.

\begin{lem}\label{Otal}Let $\lambda$ be a disc-busting measured
lamination. Then there is a cut system
$\Delta$ with respect to which $\lambda$ is in tight position, and hence $Wh(\lambda,\Delta)$ is connected and has no
cut points.
\end{lem}
\begin{proof} 
Let $\Delta$ be any cut system of $H$.
Suppose that $\lambda$ is not in tight position with respect to $\Delta$.
Then there are arcs $\alpha$ disjoint from $\lambda$ with interior disjoint from $\Delta$ and $\beta$ on a component $D$ of $\Delta$ with $\partial \alpha=\partial \beta$ such that $\alpha$ and $\beta$ are homotopic in $H$ fixing the endpoints but are not on $\partial H$.
Let $\gamma$ be $\alpha \cup \beta$.
Since $\alpha$ is homotopic to $\beta$ in $H$ fixing the endpoints, $\gamma$ bounds a compressing disc $D'$.
Since $\lambda$ was assumed to be disc-busting, $i(\lambda,\partial D') > \eta$.
On the other hand, $(\partial D \setminus \beta) \cup \alpha$ is also a meridian, and bounds a compressing disc $D''$.
Since $\beta$ is disjoint from $\lambda$, by replacing $D$ with $D''$, the intersection number with $\lambda$ is reduced by more than $\eta$.
We now define a new cut system $\Delta'$ to be $(\Delta \setminus \{D\}) \cup \{D''\}$.
As was observed above we have $i(\lambda, \Delta') < i(\lambda, \Delta)-\eta$.
Since $i(\lambda, \Delta)$ is finite and any meridian has intersection number more than $\eta$ with $\lambda$, this process must terminate in finite steps, and we reach a cut system  with respect to which $\lambda$ is in tight position.
\end{proof}

Originally in \cite{Wh1, Wh2}, Whitehead considered $Wh(A,X)$ for some finite set
$A\subset F$ and a generating set $X$ to check the separability of $A$. 
Here, $A$ is said to be
{\it separable} if there exists a free decomposition $F=F_1* F_2$ such that
for every element $a\in A$, it is conjugate into one of $F_i$.
The vertices of $Wh(A,X)$ are $\{x^{\pm}|x\in X\}$ and for any $a\in
A$, two vertices $x$ and $y$ are connected by an edge if and only if $xy^{-1}$
appears in $a$ or in a cyclic permutation of $a$. 
If $g$ is primitive, then it is separable. 
Note that for a primitive word $g\in\pi_1(H)$, the homotopy class $g$ can be realised as a simple closed geodesic $\lambda(g)$ on $\partial H$ and then $Wh(\lambda(g), \Delta)=Wh(g, X(\Delta))$ when we choose $\Delta$ so that $\lambda(g)$ does not contain any $\Delta$-wave. Here $X(\Delta)$ is the dual generating system of $F$ corresponding to $\Delta$.
Thus our previous definition for $Wh(\lambda, \Delta)$ accords with this original definition.

Due to the following lemma \cite{Wh1, Wh2}, for a primitive element $g$, $Wh(g,X)$ is either disconnected or
has a cut point for any generating set $X$.
\begin{lem}\label{Whitehead}(Whitehead) Let $g$ be
a cyclically reduced word in a free group $F$, and let $X$ be a
fixed generating set. If $Wh(g,X)$ is connected and has no cut point,
then $g$ is not separable, hence in particular $g$ is not primitive.
\end{lem}

A word $g$ is called {\it blocking}  if there exists $n$ such that $g^n$ is not a subword of a cyclically reduced primitive word. 
Likewise, a lamination $\lambda$ on $\partial H$ is called {\it blocking} with respect to a generating system if
$\lambda$ is in tight position with respect to the cut system dual to the generating system, and there exists some $k$ such
that every length $k$ subword of the infinite word determined by a
leaf of $\lambda$ does not appear as a subword in a cyclically reduced primitive
word. Recall that a  primitive word is a member of a  generating
system of $F$.

\begin{cor}\label{blocking}A disc-busting minimal lamination $\lambda$ on the
boundary of a handlebody is blocking with respect to some
generating set.
\end{cor}
\begin{proof}
Consider a generating system dual to a cut system given in Lemma \ref{Otal}.
Suppose $\lambda$ is not blocking. Then for all $k$, there is a subword of length $k$ of the infinite word determined by
$\lambda$ which is a subword of a cyclically reduced primitive word.
Taking a sufficiently large $k$, the Whitehead graph of the length $k$ subword is equal to $Wh(\lambda, \Delta)$. Lemmata \ref{tight} and \ref{Whitehead} give us a contradiction.
\end{proof}

When we consider Cannon-Thurston maps, we shall need to deal with both geodesics on $\partial H$ and those in the Cayley graph $\Gamma_F$.
For that, we shall reinterpret the lemmata above for geodesics in the Cayley graph.

Let $\widetilde H$ be the universal cover of the handlebody $H$.
For any fixed Riemannian metric on $H$, its pull back to $\widetilde H$ makes $\widetilde H$ quasi-isometric to the Cayley graph $\Gamma_F$ (for any generator system).
Therefore in particular $\widetilde H$ is Gromov hyperbolic, and it can be compactified by adding the boundary at infinity $\partial_\infty F$ of $F$.

Let $\Delta$ be a cut system for $H$.
Then the preimage $\widetilde \Delta$ of $\Delta$ in $\widetilde H$ cuts $\widetilde H$ into balls, and each of its components separates $\widetilde H$ into two.
A point in $\partial_\infty F$ corresponds to a sequence of distinct components (discs) $\{D_1, D_2,\dots \}$ in $\widetilde \Delta$ such that all of $D_{i+1}, D_{i+2}, \dots$ lie on the same side of $D_i$.
For two such sequences of discs $\{D_i\}$ and $\{D'_i\}$, they represent the same point at infinity if and only if  for each $D_i$ there is $j_0$ such that all the discs  $D'_{j_0}, D'_{j_0+1}, \dots $ lie on the same side of $D_i$ as $D_{i+1}$, and conversely for each $D'_j$ there is $i_0$ such that all the discs $D_{i_0}, D_{i_0+1}, \dots$ lie on the same side of $D'_j$ as $D'_{j+1}$.
Let us say two discs $D$ and $D'$ in $\widetilde \Delta$ are adjacent if they lie on the boundary of the same ball obtained by cutting $\widetilde H$ along $\widetilde \Delta$.
For any point at infinity $p \in \partial_\infty F$, we can choose a sequence $\{D_i\}$ as above representing $p$ such that $D_i$ and $D_{i+1}$ are adjacent for every $i$.
We call such a sequence {\it maximal}.
For two maximal sequences $\{D_i\}$ and $\{D'_j\}$, they represent the same point at infinity if and only if they are eventually the same, \ie, there are $i_0$ and $j_0$ such that $D_{i_0+k}=D'_{j_0+k}$ for every $k \in \N$. 

Let $\Gamma_F$ be a Cayley graph of $F$ with respect to the generator system dual to $\Delta$.
Since $\vee S^1$ is a spine of $H$ and we identified $\widetilde{\vee S^1}$ with $\Gamma_F$, we can regard the Cayley graph $\Gamma_F$ as being embedded in the interior of $\widetilde H$. Then there is a quasi-isometric deformation retraction $r: \widetilde H \rightarrow \Gamma_F$ which projects each disc of $\widetilde \Delta$ to a point.
Now, for any geodesic $l$ in $\Gamma_F$ which can be either a segment, a ray or a line, we can define the Whitehead graph $Wh(l, \Delta)$ in the same way as for laminations on $\partial H$ as follows.

We consider the discs in $\widetilde \Delta$ which $l$ intersects and array them in the order in which $l$ intersects them, to get a sequence $\{D_i\}$.
This sequence is maximal in the sense defined above.
Note that each disc $D_i$ in $\widetilde \Delta$ projects to a disc $\delta_i$ in $\Delta$.
We construct a Whitehead graph $Wh(l, \Delta)$ by letting the vertices be as before and connecting two vertices corresponding to $\delta_i$ to $\delta_{i+1}$ for each $i$, with signs in accordance with the way $l$ connects $D_i$ and $D_{i+1}$.

\begin{lem}
\label{lift of leaf}
Let $\lambda$ be a minimal lamination on $\partial H$ which is in tight position with respect to $\Delta$.
Let $l$ be a lift  of a leaf of $\lambda$ to $\widetilde H$, and let $l^*$ be a geodesic in $\Gamma_F$ which has the same endpoints at infinity as $r(l)$ where $r: \widetilde H \rightarrow \Gamma_F$ is the quasi-isometric deformation retraction.
Then we have $Wh(l^*, \Delta)=Wh(\lambda, \Delta)$.
\end{lem}
\begin{proof}
Since $\lambda$ is minimal, every leaf is dense in $\lambda$. Thus, if an arc in $\lambda \setminus \Delta$ connects two discs in some homotopy class of arcs (fixing the endpoints), then so does the projection of $l$ to $S$.
Therefore $l$ passes through discs of $\widetilde \Delta$ in such a way that every pair of successive intersections corresponds to an edge of $Wh(\lambda, \Delta)$ and every edge in $Wh(\lambda, \Delta)$ is realised by a pair of successive intersections. 
Therefore, we have only to show that $l^*$ and $l$ intersect the same discs in $\widetilde \Delta$ in the same order.
Since $l^*$ is a geodesic in the tree, and $l$ and $l^*$ have the same endpoints at infinity, $l$ must intersect all the discs in $\widetilde \Delta$ that $l^*$ intersects.
On the other hand $l$ can have an intersection with a disc $D$ in $\widetilde \Delta$ which is disjoint from $l^*$ only when it intersects $D$ twice in the opposite directions.
This implies that $l$ contains a wave and thus $l$ is not in tight position with respect to $\Delta$, contradicting our assumption.
Thus we have completed the proof. 
\end{proof}

We say that a geodesic ray  or a line $k$ in $\Gamma_F$ is {\em asymptotic} to a geodesic lamination $\lambda$ if there is a lift $\widetilde l$ of a leaf of $\lambda$ to $\widetilde H$ which shares at least one of the endpoints at infinity with $k$.

\begin{lem}
\label{asymptotic}
Suppose that a geodesic $k$ in $\Gamma_F$ is asymptotic to a minimal geodesic lamination $\lambda$ which is in tight position with respect to $\Delta$.
Then we have $Wh(k,\Delta)\supset Wh(\lambda,\Delta)$.
\end{lem}
\begin{proof}
Let $l$ be a lift of a leaf of $\lambda$ which shares an endpoint at infinity with $k$, and $l^*$ the geodesic in $\Gamma_F$ having the same endpoints at infinity as $l$.
By the proof of the previous lemma, we see that $l^*$ intersects the same discs in $\widetilde \Delta$ in the same order as $l$.
Let $\{D_i\}$ be a sequence of discs in $\widetilde \Delta$ arrayed in the same order as $l^*$ intersects them.
Since $k$  shares an endpoint with $l$, hence with $l^*$, there is $i_0$ such that $k$ intersects $\{D_{i_0+1},D_{i_0+2}, \dots \}$ in this order.
Now, since $\lambda$ is minimal, each leaf of $\lambda$ is recurrent.
Let $l'$ be a sub-leaf of $l^*$ starting from the intersection of $l$ with $D_{i_0}$.
Then we have $Wh(l', \Delta)=Wh(l^*,\Delta)=Wh(\lambda, \Delta)$.
Since $k$ intersects $\{D_{i_0}, D_{i_0+1}, \dots\}$ in this order, we have $Wh(k, \Delta) \supset Wh(l',\Delta)$.
This completes the proof.
\end{proof}

\subsection{Hyperbolic 3-manifold} We shall mainly concentrate on hyperbolic handlebodies which can be represented as $M=\H^3/\rho(F)$ for a discrete faithful representation $\rho$.
If $M$ is geometrically finite without parabolics, then $\rho$ is just a Schottky representation by Maskit \cite{Mas} and if $M$ is geometrically infinite without parabolics, then its compact core is a compact handlebody $H$ and the end $M\backslash int(H)$ is homeomorphic to $\partial H \times [0,\infty)$ by the Tameness theorem \cite{Ag, CG}.
In this case,  the ending lamination on $\partial H$ is connected, filling, and contained in the Masur domain of $\partial H$ by Canary \cite{Ca1}. 
When $\rho(F)$ has parabolics, they all have to be contained in  rank-1 maximal parabolic groups, and a relative compact core is also homeomorphic to a compact handlebody $H$ which meets each rank-$1$ closed cusp neighbourhood along a single annulus. The complement of these annuli in $\partial H$ consists of  several components $S_i$ which may be either compressible or incompressible. 
Each end neighbourhood $E_i=S_i\times[0, \infty)$ facing $S_i$ may be geometrically finite or otherwise be geometrically infinite, and in the latter case it has the ending lamination $\lambda_i$.
The union of a relative compact core and the closed cusp neighbourhoods is called an {\it augmented Scott core} and will be denoted by $H'$.

\section{Cannon-Thurston map}
Starting with the pioneering work of Cannon and Thurston \cite{CT} for
closed 3-manifolds fibring over a circle, 
Cannon-Thurston maps have been generalised in several ways by
Bowditch \cite{Bo}, Klarreich \cite{Kl}, McMullen \cite{Mc}, Minsky
\cite{Min1}
 and
Mj \cite{Mj1, Mj2, Mj3, Mj5, Mj6, Mj7} (see also \cite{Mi},
\cite{Sou}).
Mj proved the existence of Cannon-Thurston maps for Kleinian surface groups in \cite{Mj5}, and described the points identified by Cannon-Thurston maps  in \cite{Mj6}(see \cite{Mj9} for the case of punctured surfaces).
Recently he also proved the existence of Cannon-Thurston maps for arbitrary Kleinian groups and described the points identified by the maps in \cite{Mj8}. 
We shall make use of his result in the case of free Kleinian groups as a main ingredient for the proof of our main results.

Given a discrete faithful representation $\rho : F \rightarrow\pslc$ without parabolics, the {\it Cannon-Thurston map} is a continuous extension of
$\tau_{\rho,o}:\Gamma_F \ra \H^3$ to  $\widehat\tau_{\rho,o}:\widehat\Gamma_F\ra \H^3\cup\hat\C$
where $\hat\C$ is the ideal boundary of $\H^3$ and $\widehat\Gamma_F$ is the Gromov compactification of $\Gamma_F$. 
The boundary of $\widetilde H$ is a covering of $S$, which we denote by $S_F$.
We fix a hyperbolic metric on $S$ and pull it back to $S_F$.
As was explained in \S \ref{sec: Whitehead}, the Gromov boundary of both $\widetilde H$ and $\widehat\Gamma$ is identified with $\partial_\infty F$.
Mj's theorem for free Kleinian groups without parabolics is as follows.
\begin{thm}\label{Mj1}
Given a discrete faithful representation $\rho:F\ra \pslc$ without parabolics, let $\Gamma_F$ be the Cayley graph of $F$ and $i:\Gamma_F\ra \H^3$  the natural identification of $\Gamma_F$ with its image under $\tau_{\rho,o}$ for a chosen base point $o\in\H^3$. 
Let $\lambda$ be the ending lamination of $\rho(F)$ and $\widetilde \lambda$ its preimage in $S_F$.
Then $i$ extends continuously to a map $\widehat{i}:\widehat{\Gamma}_F\ra \H \cup \hat \C$.
If we let $\partial i$ denote the restriction of $i$ to the Gromov boundary $\partial_\infty F$, then $\partial i(a)=\partial i(b)$ if and only if $a, b$ are either ideal endpoints of a leaf of $\widetilde \lambda$, or ideal vertices of one of the complementary ideal polygons of $\widetilde \lambda$, where $\partial_\infty F$ is regarded as the boundary at infinity of $S_F$.
\end{thm}

Recall that $M:=\H^3/\rho(F)$ can be decomposed into $M=H\cup E$ where $H$ is a genus $n$ compact handlebody and $E$ is homeomorphic to $\partial H\times[0, \infty)$. We use the convention that $\partial H$ and $\partial H\times\{0\}$ are identified and denote $\partial H$ as $S$. 
When the end $E$ is geometrically infinite, we have an ending lamination $\lambda$ on $S$. By \cite{Ca1}, $\lambda$ is in Masur domain and by \cite[Theorem 1.3]{Ot}, there is a cut system $\Delta$ with respect to which $\lambda$ is in tight position (the arguments of the proof of \cite[Theorem 1.3]{Ot} are the arguments we used in the proof of Lemma \ref{Otal}).
The lamination is realised as a geodesic lamination uniquely once we fix a hyperbolic metric on $S$.
As was explained in \S \ref{sec: Whitehead}, each leaf of $\lambda$ can be lifted to a geodesic $l$ in $\partial \widetilde H$ (with the pulled-back hyperbolic metric) and its projection by the retraction $r$ forms a quasi-geodesic in $\Gamma_F$.
We denote the geodesic on $\Gamma_F$ with the same endpoints at infinity as $r(l)$ by $l^*$ as before.
Then $Wh(l^*, \Delta)=Wh(\lambda, \Delta)$ by Lemma \ref{lift of leaf}.
By Lemma \ref{asymptotic}, for any geodesic in $\Gamma_F$ sharing an endpoint with $l^*$, we have $Wh(k, \Delta)\supset Wh(\lambda, \Delta)$.
If we connect two endpoints $a, b$ in $\partial_\infty F$ such that $\partial i(a)=\partial i(b)$, then by Theorem \ref{Mj1}, the geodesic on $\Gamma_F$ connecting $a$ with $b$ is asymptotic to a leaf of $\widetilde \lambda$.
Therefore, we get the following corollary.

\begin{cor}\label{lem}
Given a discrete faithful representation $\rho:F\ra \pslc$ without parabolics, denote by $\lambda$ its ending lamination.
Let $\Delta$ be a cut system with respect to which $\lambda$ is in tight position.
Given $a,b\in \partial_\infty F$ which are identified by $\partial i$, if we let the  geodesic on $\Gamma_F$ joining $a,b$ be $k$, then $Wh(k,\Delta)$ contains $Wh(\lambda,\Delta)$, is connected and has no cut point. 
\end{cor}

Now we discuss the case of a free Kleinian group with parabolics. 
Recall that when $\rho(F)$ has parabolics, $M=\H^3/{\rho(F)}$ has a relative compact core $H$ of the non-cuspidal part $M_0$ which intersects each cusp neighbourhood along an annulus whose core curve we call a {\em parabolic curve}.
Each frontier component $S_i$ of $H$ in $M_0$ faces an end neighbourhood $E_i$. If $E_i$ is geometrically infinite, it has an ending lamination $\lambda_i$ which is a filling minimal lamination on $S_i$.
Now we consider the union $\Lambda$ of all the ending laminations $\lambda_i$ and all the parabolic curves.
Then $\Lambda$ itself is a  geodesic lamination on $S$.
Suppose, as in the setting of Theorem \ref{main with parabolics}, that all the parabolic curves  and the $\lambda_i$ are disc-busting.
Then $\Lambda$ is also disc-busting.
By Lemma \ref{tight}, we can find a cut system $\Delta$ with respect to which $\Lambda$ is in tight position. 
We lift $\Lambda$ to a geodesic lamination $\widetilde \Lambda$ on $S_F$. 
Since $\Lambda$ is in tight position with respect to $\Delta$, each $\lambda_i$ has no $\partial\Delta$-waves although $\lambda_i$ may not be in tight position with respect to $\Delta$. Thus no leaf of $\widetilde \Lambda$ can intersect a component of the preimage of $\Delta$ twice and each leaf of $\widetilde \Lambda$ connects two distinct points at infinity on $\partial_\infty F$.

We define an equivalence relation $\widetilde{\mathcal R}$ on $\partial_\infty F$ by letting $a \widetilde{\mathcal R} b$ if and only if one of the following holds :
\begin{enumerate}
\item $a$ and $b$ are the endpoints of a leaf of $\widetilde\Lambda$;
\item there is a complementary region of $\widetilde \Lambda$ with trivial or cyclic stabiliser whose closure contains both $a$ and $b$;
\item there are two complementary regions $U,U'$ of $\widetilde \Lambda$ whose closures share a lift of a parabolic curve, such that $a$ is contained in the closure of $U$ and $b$ is contained in the closure of $U'$.
\end{enumerate}

Here, the stabiliser refers to that with respect to the action of $F$ on $S_F$. 
A complementary region of $\widetilde \Lambda$ which is a lift of a complimentary region of $\Lambda$ bounded only by parabolic curves may have non-cyclic stabiliser.
Such a region corresponds to a geometrically finite end of $M_0$

It is easy to check that this relation $\widetilde{\mathcal R}$ is the transitive closure of another relation ${\mathcal R}$ which is defined by letting $a \mathcal R b$ if and only if $a$ and $b$ are either the endpoints of a leaf of $\widetilde\Lambda$ or two vertices of the same complementary region of $\widetilde \Lambda$ with trivial or cyclic stabiliser. 

Mj's Theorem about the identified points of a Cannon-Thurston map can be adapted to our case as follows. Recall that $H'$ is an augmented Scott core of $M$. 
\begin{thm}\label{Mj2}
Let $\rho:F\ra\pslc$ be a discrete faithful representation. 
Then the natural identification $i$ of $\Gamma_F$ with its image under $\tau_{\rho,o}$ extends continuously to a map $\widehat i : \Gamma_F \cup \partial_\infty F \rightarrow \H^3 \cup \hat \C$.
Furthermore for $a,b\in\partial_\infty F$, $\partial i(a)=\partial i(b)$ if and only if $a\widetilde{\mathcal R} b$.  
\end{thm}

Assume that $\partial i(a)=\partial i(b)$ for $a,b\in \partial_\infty F$. 
Let $k$ be a geodesic connecting $a$ and $b$ on $\Gamma_F$.
Then $k$ is asymptotic to some leaf $l$ of $\widetilde \Lambda$ by our definition of $\mathcal R$.
Let $\lambda$ be a component of $\Lambda$ containing the projection of $l$ and take $\Delta'$ to be the disc system such that $Wh(\lambda, \Delta')$ is connected and has no cut point coming from Lemma \ref{tight}. Since two Cayley graphs coming from two disc systems $\Delta'$ and $\Delta$ are quasi-isometric, if we let $k'$ be the geodesic in the Cayley graph with respect to $\Delta'$ connecting $a$ and $b$, then $k'$ is asymptotic to $l$.  
By the same argument as in the case without parabolics, Corollary \ref{lem} holds also in this case.

\begin{cor}\label{lem with parabolics}
Given a discrete faithful representation $\rho:F\ra \pslc$, denote by $\Lambda$ the union of its ending laminations and parabolics. Suppose that all the components of $\Lambda$ are disc-busting.
Given $a,b\in \partial_\infty F$ which are identified by $\partial i$, if we let the  geodesic on $\Gamma_F$ joining $a,b$ be $k$,  then $Wh(k,\Delta)$ contains $Wh(\lambda,\Delta)$ for some component $\lambda$ of $\Lambda$ and a cut system $\Delta$ with respect to which $\lambda$ is in tight position. Then $Wh(k,\Delta)$ is connected and has no cut point for the disc system $\Delta$.  
\end{cor}

\section{Primitive stable representations of a free group}
\subsection{Free groups without parabolics}
In this section, we shall prove the first of our theorems.
The overall argument is as follows. 
Assuming that the representation is
not primitive stable, there exists a sequence of primitive
elements $\{w_n\}$ such that $\{\tau_{\rho,o}(\widetilde w_n)\}$ is not a family of
uniform quasi-geodesics, where $\widetilde w_n$ is a line in $\Gamma_F$ corresponding to $w_n$. 
After passing to a subsequence, $\widetilde w_n$ converges uniformly on every compact set to a bi-infinite line $\widetilde w_\infty$ in the
Cayley graph with two distinct end points in the Gromov boundary.
Then we shall show that the endpoints of $\tau_{\rho,o}(\widetilde w_\infty)$
are the same point in $\partial \H^3=\hat \C$.
Applying Mj's result cited as Theorem  \ref{Mj1} above, we shall conclude that the
endpoints of $\widetilde w_\infty$ are either the endpoints of a lift of a leaf of
the ending lamination of $M=\H^3/\rho(F)$ or ideal end-points of a complementary polygon. 
Finally we shall apply Lemma \ref{Whitehead} to draw a contradiction.

\begin{proof}[Proof of Theorem \ref{main without parabolics}]
Suppose that $\rho$ is not primitive stable.
Then there exists a sequence $\{w_n\}$ of cyclically reduced primitive words such that the $\tau_{\rho,o}(\widetilde w_n)$  are not uniform quasi-geodesics for every lift $\widetilde w_n$ of $w_n$ to $\Gamma_F$. 
Recall that $M=\H^3/\rho(F)$ is homeomorphic to $H\cup (\partial H\times[0,\infty))$ where $\partial H$ is identified with $\partial H\times\{0\}$. 

Let $\gamma_{w_n}$ be the geodesic in $\H^3$ joining the endpoints at infinity of $\tau_{\rho,o}(\widetilde w_n)$. 
Note that $\gamma_{w_n}$ is the axis of the loxodromic isometry $\rho(w_n)$.
By the following Lemma \ref{4.1} which was observed by Minsky in \cite[Lemma 3.2]{Min2}, it is not possible that all $\gamma_{w_n}$ are contained in a uniformly bounded
neighbourhood of $\tau_{\rho,o}(\Gamma_F)$.

\begin{lem}\label{4.1}
A discrete faithful representation $\rho:F\rightarrow \pslc$ is primitive stable if and only if there is uniform neighbourhood of $\tau_{\rho,o}(\Gamma_F)\subset M=\H^3/\rho(F)$ containing all the closed geodesics representing primitive elements.
\end{lem}
\begin{proof}
If $\rho$ is primitive stable, there are $K, \epsilon$ such that for any primitive word $w$, the line $\tau_{\rho,o}(\widetilde w)$ is a $(K,\epsilon)$-quasi-geodesic. This implies that there is a uniform $L$ such that the line $\tau_{\rho,o}(\widetilde w)$ and the geodesic line connecting its ideal endpoints are within Hausdorff distance $L$. By projecting them down to $M$, we see that any primitive closed geodesic stays in the $L$-neighbourhood of the projection of $\tau_{\rho,o}(\Gamma_F)$.

To prove  the converse, we can assume the projection to $M$ of the uniform neighbourhood of $\tau_{\rho,o}(\Gamma_F)$ to be a core handlebody of $M$, and regard it as $H$.
Let $r : H \ra \vee S^1$ be the obvious retraction and let $\widetilde r : \widetilde H \ra \Gamma_F$ be its lift.
Then $r$ and $\widetilde r$ are $(K,\delta)$ quasi-isometries for some constants $K,\delta$ depending only on $H$ and $r$. Thus we get a $(K,\delta)$ quasi-geodesic $\widetilde{r}(\gamma_{w})$. Note that $\widetilde{r}(\gamma_{w})=\tau_{\rho,o}(\widetilde w)$ because their projections to $M$ are freely homotopic in $H$.
\end{proof}

Hence, there exists a sequence of positive numbers $\left\{\epsilon_n\right\}$ such that $\gamma_{w_n}$ is
not contained in the $\epsilon_n$-neighbourhood of the core graph $\tau_{\rho,o}(\Gamma_F)$
where $\epsilon_n\rightarrow \infty$. 
Since $\gamma_{w_n}$ is not
contained in the $\epsilon_n$-neighbourhood of
$\tau_{\rho,o}(\Gamma_F)$, neither is it in
the $\epsilon_n$-neighbourhood of $\tau_{\rho,o}(\widetilde w_n)$.
Thus, we can choose $p_n\in\gamma_{w_n}$ which is not contained in the $\epsilon_n$-neighbourhood of 
$\tau_{\rho,o}(\widetilde w_n)$. Then the geodesic plane orthogonal to $\gamma_{w_n}$ at $p_n$ has to intersect $\tau_{\rho,o}(\widetilde w_n)$ on at least one point, say $q_n$, so that we get $d_{\H^3}(q_n, \gamma_{w_n})>\epsilon_n$.
If $q_n$ lies in a geodesic segment $s$ of $\tau_{\rho,o}(\widetilde w_n)$, then one of the two endpoints of $s$, which we denote by $\mathbf v$, has distance from $\gamma_{w_n}$ larger than $\epsilon_n$ because the minimal distance function $d_{\H^3}(q, \gamma_{w_n})$ where $q$ moves in $s$ is convex (see \cite[pp. 178]{BH}, for instance).
Moreover, we can shift the words $w_n$ so that the specified vertex $\mathbf v$ becomes the base point $o$ as
follows.

The vertex $\mathbf v$ on $\tau_{\rho,o}(\widetilde w_n)$ is expressed as $\rho(w_n^i v_n)o$ where $w_n=g_1g_2\ldots g_k$
and $v_n=g_1\ldots g_l$ for $l<k$ and $i\in \Z$. 
Assume that $d_{\H^3}(\rho(w_n^i v_n)o,\gamma_{w_n}) > \epsilon_n$.
Then, since $\gamma_{w_n}$ is the axis of the loxodromic isometry
$\rho(w_n)$, we get
$$d_{\H^3}(\rho(w_n^i v_n)o,\gamma_{w_n})=d_{\H^3}(\rho(v_n)o,\gamma_{w_n})=d_{\H^3}(o,\rho(v_n)^{-1}\gamma_{w_n}),$$
and
$$\rho(v_n)^{-1}{\gamma_{w_n}}=\gamma_{v_n^{-1}w_nv_n}.$$ 
Since 
$v_n^{-1}w_nv_n$ is a shifted word of $w_n$, it is also cyclically reduced
and primitive. 
We also have
$$d_{\H^3}(o,\gamma_{v_n^{-1}w_nv_n}) > \epsilon_n.$$
Therefore, $\gamma_{v_n^{-1}w_nv_n}$ has to leave every compact subset in
$\H^3$ as $n \rightarrow \infty$, and $\{v_n^{-1}w_nv_n\}$ is a desired sequence.

Now, we have obtained a new sequence $\left\{w_n'\right\}$ such that
$d_{\H^3}(o, \gamma_{w_n'})$ goes to $\infty$ as $n\rightarrow
\infty$, which means that the spherical distance between the endpoints of $\gamma_{w_n'}$ goes to $0$,  so that the two endpoints of
$\tau_{\rho,o}(\widetilde w_n')$ in $\partial{\H^3}$ converge to
the same point as $n \rightarrow \infty$.

\begin{lem}
\label{Fl}
In $\Gamma_F$, after passing to
a subsequence, $\{\widetilde w_n'\}$ converges uniformly on every compact set to a bi-infinite geodesic
$\widetilde w_\infty$ with distinct end points and such that for some cut system $\Delta$,
we have $Wh(\widetilde w_\infty,\Delta)\supset Wh(\lambda,\Delta)$, where $\lambda$ is the ending lamination of $\rho(F)$. 
\end{lem}
\begin{proof}Since $F$ has a set of finite number of generators $X$, after passing to a
subsequence, we may assume that all $w_n'$ have $x_i$ and $x_j$ in $\{x^{\pm}|x\in X\}$ as their first
and last letter respectively for fixed $i,j$. Note that $x_j$ cannot be $x_i^{-1}$ because $w'_n$ is cyclically reduced. 
Then the geodesic $\widetilde w_n'$ has endpoints at infinity
in the regions determined by $x_i$ and $x_j^{-1}$.
Hence, passing to a subsequence, $\widetilde w_n'$ has a limit geodesic $\widetilde w_\infty$
with distinct endpoints $a,b$ at infinity in the regions
determined by $x_i$ and $x_j^{-1}$. 
Since the endpoints of
$\tau_{\rho,o}(\widetilde w_n')$ in $\partial{\H^3}$ converge to
the same point, $a$ and $b$ are identified under the Cannon-Thurston
map. 

\end{proof}
Now returning to the main proof,
$Wh(\lambda,\Delta)$ is connected and has no cut point with respect
to  $\Delta$ by Lemma \ref{lem} and from $Wh(w_\infty,\Delta)\supset Wh(\lambda,\Delta)$, we can see that the same is true for
$Wh(w_n',\Delta)$ for large $n$. 
On the other hand for any primitive word $w_n'$,
this is not possible by Corollary \ref{Whitehead}.
\end{proof}

\subsection{Free groups with parabolics}
Recall that when $\rho:F\ra \pslc$ has parabolics,  $M=\H^3/\rho(F)=H'\cup\cup E_i$ where $H'$ is the augmented scott core and
$E_i$ is an end neighbourhood facing the relative compact core $H$ along $S_i\subset \partial H$. 
The $S_i$  are glued to each other by parabolic loci $A_i$  whose core curves are denoted by $c_i$. 
When $E_i$ is geometrically infinite, it has ending lamination $\lambda_i$ which is disc-busting by assumption.

Before we start the proof of our second theorem, we remark that
in the case when $\rho$ has parabolics,  the
manifold may not be primitive stable even if every proper free factor is Schottky. 
The following example due to Minsky shows this. 
Let $M$ be a handlebody of even genus. Then $M$ is homeomorphic to $\Sigma\times I$ where $\Sigma$ is a
genus $g$ surface with one boundary component. Consider a discrete faithful representation $\rho: \pi_1(M)=F_{2g}\rightarrow \pslc$ such that the
boundary curve of $\Sigma$ corresponds to a parabolic element and at least
one end is degenerate. Then we can see that $\rho$ is not
primitive stable using the following argument which is also due to Minsky.

First, we can see that the covering of $M$ corresponding to a free factor of $F_{2g}$ is convex cocompact by using Canary's covering theorem \cite{Ca2}.
This shows that if we restrict $\rho$ to a proper free factor, then the representation is Schottky.
Noting that every non-peripheral non-separating simple closed
curve on $\Sigma$ is primitive,  
suppose that $\{p_i\}$ is a sequence of such primitive
simple closed curves converging to the ending lamination of $M$, whose geodesic representatives exit  the end. If we suppose $\rho$ is primitive stable, 
then a line passing through the identity
element determined by the conjugacy class of $p_i$ in the Cayley
graph is mapped to a uniform $(K,d)$-geodesic in $\H^3$ passing
through a fixed point $o$ since the $p_i$ are primitive. Then the geodesic lines homotopic to its image by $\tau_{\rho,o}$ cannot be
far from $o$, hence their images in $M$ are near the projection of
$o$.
This contradicts the fact that they exit the end.

We note that in this example the ending lamination, which we denote by $\lambda$, is not disc-busting.
In fact, $\lambda$ is contained in the Hausdorff limit of meridians of the form of $a_i \times \{0\} \cup \partial a_i \times I \cup a_i \times \{1\}$, where the $a_i$ are essential arcs on $\Sigma$ whose Hausdorff limit contains $\lambda$.

\begin{proof}[Proof of Theorem \ref{main with parabolics} (sufficiency)]
If there is no ending lamination $\lambda_i$, then by the same argument as the proof of Theorem 4.1 in  Minsky \cite{Min2}, we see that $\rho$ is primitive stable. 

If there is at least one ending lamination, then we consider $\Lambda$ which is the union of the ending laminations and parabolic curves, and get the equivalence relation $\widetilde{\mathcal R}$ from $\widetilde \Lambda$ as we explained just before stating Theorem \ref{Mj2}.  Repeating the same argument as in the case without parabolics replacing Corollary \ref{lem} with Corollary \ref{lem with parabolics}, we complete the proof.
\end{proof}

By Minsky's result that the restriction of a primitive stable
representation to a proper free factor of the free group is
Schottky, we immediately get the following corollary.
\begin{cor}Let $M=\H^3/\rho(F)$ as in Theorem \ref{main without parabolics} or \ref{main with parabolics}.
 If $F_n=A*B$ into two nontrivial free factors, then the
covering manifold corresponding to $A$ or $B$ is convex cocompact,
i.e., Schottky.
\end{cor}
Note that the the above corollary can be also obtained by Canary's covering theorem (\cite{Ca2}) directly.

\section{Necessary condition}
In this section, we shall prove the necessity part of Theorem \ref{main with parabolics}. Actually we shall prove the contrapositive, namely we assume that some parabolic or ending lamination is disc-dodging and show that the representation is not primitive stable. The case that $\rho$ has a parabolic curve is dealt with by the next Lemma \ref{parabolic locus}.

\begin{lem}\label{parabolic locus}
If there is a parabolic curve on $\partial H$ which is disjoint from a meridian, then $\rho$ is not primitive stable.
\end{lem}

\begin{proof}
Suppose that $c$ is a parabolic curve which is disjoint from a meridian $m$.
If $m$ is separating, then it bounds a separating compressing disc $D$ which gives rise to a free-product decomposition of $F=\pi_1(H)$.
One of the free factors which contains an element corresponding to $c$ is not Schottky.
By Minsky's result in \cite{Min2} stated as (3) in our Introduction, this shows that $\rho$ is not primitive stable.

In the case when $m$ is non-separating,  we have a decomposition of $F$ of the form $F=A*\Z$ such that $c$ represents an element conjugate into $A$. The same result of Minsky leads to the same conclusion, i.e. $\rho$ is not primitive stable.
\end{proof}

To conclude the proof of the necessity part, we need to assume that an ending lamination $\lambda$ is disc-dodging and conclude that the representation is not primitive stable. In order to do that, we shall approximate $\lambda$ by disc-dodging curves. The following Lemma \ref{esan}, which is a standard result of $3$-manifold topology, will help us construct disc-dodging curves.

\begin{lem}		\label{esan}
For any essential annulus $A$ in $H$, there is a meridian which is disjoint from $A$.
\end{lem}

\begin{proof}
What we want to show is that there is a compressing disc for $H$ disjoint from $A$.
Let $D$ be a compressing disc for $H$.
We isotope $D$ so that there is no inessential intersection between $D$ and $A$.
If $D \cap A= \emptyset$, then we are done.
Suppose not.
We consider an arc $k$ in $D \cap A$ which is outermost in $D$ and cuts off a semi-disc $\Delta$ from $D$.
If $k$ connects the same component of $\partial A$, then it cuts off a disc $\Delta'$ from $A$ and $\Delta \cup \Delta'$ is a compressing disc which can be isotoped off $A$.
If $k$ connects two components of $\partial A$, then we can boundary-compress $A$ along $\Delta$, and get a compressing disc disjoint from $A$.
%
\end{proof}

As was mentioned before, we shall approximate a disc-dodging minimal lamination by disc-dodging curves, but before going any further, we need a few extra definitions.

We fix some hyperbolic structure on $\partial H$.
For a minimal lamination $\mu$ which is not a simple closed curve we denote  by $S(\mu)$ the unique minimal compact subsurface of $\partial H$ with geodesic boundaries containing $\mu$.
In $S(\mu)$ there is a unique maximal multi-curve that is disjoint from $\mu$, we denote it by $C_\mu$. We have $\partial S(\mu)\subset C_\mu$ but this inclusion may not be an equality.


\begin{lem}
\label{ending lamination}
Let $\lambda$ be a disc-dodging minimal lamination which is not a simple closed curve, then either every component of $C_\lambda$ is disc-dodging or $\lambda$ is contained in  a Hausdorff limit of a sequence of disc-dodging simple closed curves on $\partial H$ as its unique minimal component.
\end{lem}

\begin{proof}
The idea of the proof is to use a homoclinic leaf (see the definition below) to construct a sequence of essential annuli whose boundaries approximate $\lambda$ as in \cite[Lemme C.1]{Le} (see also the proof of \cite[Proposition 1]{KS}) and to use Lemma \ref{esan}. Although we cannot adapt the arguments of \cite[Lemme C.1]{Le} to our case, it will turn out that we can modify our settings so that we can apply some of the results of \cite{Le}.

Consider a disc-dodging minimal lamination $\lambda$ endowed with a transverse measure and a sequence of meridians $m_i$ such that $ i(m_i,\lambda)\ra 0$. If all the components of $C_\lambda$ are disc-dodging, we are done. Thus, we assume that there is a disc-busting component $c$ of $C_\lambda$. 
We endow each leaf of $C_\lambda$ with a transverse Dirac mass with a weight equal to $\pi$ and denote by $\gamma$ a measured lamination which is the union of the weighted multi-curve thus obtained and $\lambda$. 

If an essential annulus $A$ satisfies $i(\partial A,\gamma)=0$, we cut $H$ along $A$ and keep the component  which contains $\lambda$, and denote it by $H_1$. 
Notice that $H_1$ is also a handlebody. Applying the cut-and-paste operation described in the proof of Lemma \ref{esan} to $A$ and a disc bounded by $m_i$, we get a meridian $m_{i,1}\subset\partial H$ which is disjoint from $A$ and satisfies $i(\lambda,m_{i,1})\leq 2 i(\lambda, m_i)$. 
Since $c\subset C_\lambda$ is disc-busting, $m_{i,1}$ intersects $c$. In particular, $m_{i,1}$ lies on $\partial H_1$, which implies that $\lambda$ is disc-dodging also in $H_1$. 
The boundary of $H_1$ contains one or two annuli corresponding to $A$ (depending on whether or not $A$ separates $H$). Inside each of these annuli, we put its core  curve on which we give  a transverse Dirac mass with  weight equal to $\pi$, and add it to $\gamma$. 
Since $m_{i,1}$ is disjoint from $A$, it is disjoint from the added components of $\gamma$.

If there is another essential annulus $A_1\subset H_1$ satisfying $i(\partial A_1,\gamma)=0$, then it can be regarded as an essential annulus in $H$ because we extended $\gamma$ so that it contains the core curve of $A$. If we cut $H_1$ along $A_1$, take the component $H_2$ containing $\lambda$, and extend $\gamma$ as before, then we get a meridian $m_{i,2}\subset \partial H_2$ with $i(\lambda, m_{i,2})\leq 2i(\lambda, m_{i,1})$. 

We repeat this operation as long as there is an essential annulus not intersecting $\gamma$ on the obtained handlebody. 
Since the number of non-parallel essential annuli in $H$ is finite, after finitely many steps, we get a new handlebody $H_\infty$ and a measured geodesic lamination $\gamma$ such that any essential annulus $A^*$ in $H_\infty$ satisfies $i(\partial A^*,\gamma)>0$. 
Furthermore there is a sequence of meridian $m_{i,\infty}\subset\partial H_\infty$ such that $i(m_{i,\infty},\lambda)\rightarrow 0$ and  $i(m_{i,\infty},\gamma)=i(m_{i,\infty},\gamma\cap S(\lambda))$.
By an abuse of notation, we denote $m_{i, \infty}$ again by $m_i$.

To apply \cite[Proposition 3.5]{Le} to our $\gamma$, we verify that $\gamma$ satisfies its hypotheses.

\begin{claim}		\label{a and c}
The measured geodesic lamination $\gamma$ satisfies :
\begin{enumerate}[(a)]
\item  the weight of any closed leaf of $\gamma$ is at most $\pi$;
\label{a}
\setcounter{enumi}{2}
\item  $i(\gamma,\partial D)> 2\pi$ for any compressing disc $D$ for $H_\infty$.
\label{c}
\end{enumerate}

\end{claim} 

\begin{proof}
The condition (\ref{a}) follows directly from the construction of $\gamma$.

By assumption a component $c$ of $C_\lambda$ is disc-busting. By definition any simple closed curve that intersects $c$ also intersects $\lambda$.  Since $c$ carries a weight equal to $\pi$ in $\gamma$ any simple closed curve $d$ that intersects $c$ at least twice satisfies $i(\gamma,d)> 2\pi$. Thus to prove that $\gamma$ satisfies the condition (\ref{c}), we need to prove that any meridian intersects $c$ at least twice.

Consider a compressing disc $D$ for $H_\infty$. Since $c$ is disc busting, $\partial D$ intersects $C$ at least once. If $\partial D$ intersects $c$ only once, we denote by $V$ a small regular neighbourhood of $c \cup D$, which is a solid torus. 
Let $D'$ be the closure of $\partial{V}\setminus \partial H_\infty$.
Then $D'$ is a compressing disc for $H_\infty$ which does not intersect $c$ since $H_\infty$ itself is not a solid torus. 
This contradicts the disc-busting property of $c$.
Thus we have proved that any meridian intersects $c$ at least twice.
\end{proof}

It follows then from \cite[Proposition 3.5]{Le} that the following two conditions are equivalent.
\begin{enumerate}[($b_1$)]
\item  There exists a positive number $\eta$ such that $i(\gamma,\partial A)\geq\eta$ for any essential annulus $A$ in $H_\infty$.
\label{b_1}
\item Let $l^+,l^-\subset\partial H_\infty\setminus \gamma$ be two disjoint half geodesics such that some lifts $\widetilde l^+$ and $\widetilde l^-$ of $l^+$ and $l^-$ to $\widetilde H_\infty$ have the same endpoints.
Then $\widetilde l^+$ and $\widetilde l^-$ are asymptotic on $\partial\widetilde H_\infty$.
\label{b_2}
\end{enumerate}

Recall that we have a sequence of meridian $m_i\subset\partial H_\infty$ such that $i(m_i,\lambda)\rightarrow 0$ and that $i(m_i,\gamma)=i(m_i,\gamma\cap S(\lambda))$. Take a subsequence so that  $\{m_i\}$ converges in the Hausdorff topology to a geodesic lamination $\mu$ on $\partial H_\infty$. By Casson's criterion (see Casson-Long \cite{CL}, Otal \cite{Ot}, and \cite[Theorem B1]{Le}), $\mu$ contains a homoclinic leaf $h$. A homoclinic leaf is defined as follows. We fix a Riemmanian metric on $H_\infty$ and lift it to $\widetilde H_\infty$. A leaf $l$ is said to be {\it homoclinic} if there are sequence of points $\{x_i\}, \{y_i\}$ on a lift $\widetilde l\subset\widetilde H_\infty$ of $l$ such that the distance between $x_i$ and $y_i$ (in $\widetilde H_\infty$) is bounded whereas their distance on $\widetilde l$ goes to $\infty$ as $i \rightarrow \infty$. 

If $h$ intersects $\lambda$ transversely, then it contains an arc $\kappa$ with  $\int_\kappa d\lambda >0$. Since $h$ lies in the Hausdorff limit of $m_i$, for $i$ large enough, $m_i$ contains an arc very close to $\kappa$. It follows that $i(m_i, \lambda)\geq \int_\kappa d\lambda$ for large $i$, where the righthand is a positive constant independent of $i$. Since $i(m_i,\lambda)\ra 0$, this cannot happen, and hence $h$ does not intersect $\lambda$ transversely.

Since $i(m_i,\gamma)=i(m_i,\gamma\cap S(\lambda))$, using the same argument, we see that $h$ is disjoint from the leaves of $\gamma$ lying outside $S(\lambda)$.  It follows that $h$ contains two disjoint half-leaves $h^+$ and $h^-$ which are disjoint from $\gamma$. 
Using these half-leaves $h^+$ and $h^-$, we shall show that $\gamma$ does not satisfy the condition $(b_2)$.  

Let $\widetilde h\subset\widetilde H_\infty$ be a lift of $h$ and $\widetilde h^\pm\subset\widetilde h$ lifts of $h^\pm$ respectively. Since $c$ is disc-busting, $\partial H_\infty \setminus c$ is incompressible. A lift of a component of $\partial H_\infty \setminus c$ to $\partial \widetilde H$ is an open disc whose closure in $\partial\widetilde H_\infty \cup\partial_\infty F'$ is a closed disc (see \cite[Lemme 2.4]{Le}), where $F'$ is a subgroup of $F$ corresponding to $\pi_1(H_\infty)$. It follows that $\widetilde h^\pm$ has a well-defined endpoint in $\partial_\infty F'$. Since $h$ is homoclinic, $\widetilde h^+$ and $\widetilde h^-$ have the same endpoint and we have the following.

\begin{claim}		\label{nota}
The half-geodesics $\widetilde h^+$ and $\widetilde h^-$ are not asymptotic on $\partial\widetilde H_\infty$
\end{claim}

\begin{proof}
Seeking a contradiction, suppose that $\widetilde h^+$ and $\widetilde h^-$ are asymptotic on $\partial\widetilde H_\infty$.
Then there is a sequence of geodesic arcs $\widetilde k_n\subset\partial\widetilde H_\infty$ joining $\widetilde h^+$ to $\widetilde h^-$ such that $\widetilde k_n\cap\widetilde h=\partial\widetilde k_n$ and  the length of $\widetilde k_n$  goes to $0$ (with respect to the pull-back of some fixed hyperbolic metric on $\partial H_\infty$). 
Let $k_n$ be the projection of $\widetilde k_n$ to $\partial H_\infty$.
Then we see that $\int_{k_n} d\lambda\rightarrow 0$, where $d \lambda$ denotes the transverse measure of $\lambda$, as follows. 
Suppose, seeking a contradiction, that there exists $\epsilon>0$ such that $\int_{k_n} d\lambda > \epsilon$.
Since $\partial H_\infty$ is compact, passing to a subsequence, the arcs $k_n$  converge to a point $p$ on $\partial H_\infty$,
and any transverse arc passing through $p$ has to have measure greater than $\epsilon$.
This contradicts the fact that $\lambda$ is minimal and is not a simple closed curve.

The arcs $k_n$ can be assumed to be homotopic without their endpoints passing through $\lambda$ since $\widetilde h^+$ and $\widetilde h^-$ are asymptotic on $\partial \widetilde H_\infty$. 
Hence, $\int_{k_n} d\lambda$ does not depend on $n$, and it follows that $\int_{k_n} d\lambda=0$ for every $n$, which means that  $k_n$ is disjoint from $\lambda$ for every $n$. 
By shortening $h^\pm$ if necessary, we may assume that $\partial\widetilde h^\pm$ lies on $\partial\widetilde k_1$. 
Then $\widetilde d=\widetilde k_1\cup (\widetilde h \setminus \widetilde h^\pm)$ bounds a disc in $\widetilde H_\infty$ which is disjoint from the preimage of $\lambda$  since $\widetilde H_\infty$ is simply connected. 
Since $k_1$ can be homotoped to arbitrarily short geodesic arcs keeping its endpoints on $\widetilde h$, it cannot be homotopic to $\widetilde h \setminus \widetilde h^\pm$, which means that $\widetilde d$ is essential on $\partial \widetilde H_\infty$.
The projection $d$ of $\widetilde d$ to $\partial H_\infty$ bounds a (possibly not embedded) disc which is disjoint from $\lambda$. 
It follows then from Dehn's lemma that there is a compressing disc for $H_\infty$ which is disjoint from $\lambda$. 
Hence each component of $C_\lambda$ is disc-dodging. This contradicts the assumption that $c\subset C_\lambda$ is disc-busting.
\end{proof}

It follows from Claim \ref{nota} that $\gamma$ does not satisfy condition $(b_2)$ either. By \cite[Proposition 3.5]{Le} there is a sequence of essential annuli $A_i$ in $H_\infty$ such that $i(\gamma,\partial A_i)\rightarrow 0$. By construction, we have $i(\gamma,\partial A)>0$ for any essential annulus $A$ (this is the cut-and-paste operation achieved for $H_\infty$). Since each leaf of $\gamma$ that is not a leaf of $\lambda$ has a weight equal to $\pi$, for $i$ large enough, $\partial A_i$ intersects $\lambda$ and $i(\partial A_i,\gamma\setminus\lambda)=0$. Extract a subsequence such that $\partial A_i$ converges in the Hausdorff topology to a geodesic lamination $\nu$. Since $i(\gamma,\partial A_i)\rightarrow 0$, the lamination $\nu$ does not intersect $\gamma$ transversely. Since $\partial A_i$ intersects $\lambda$, the only possibility is that $\lambda$ is a sublamination of $\nu$.
This implies that $\lambda$ and $\nu$ coincide as geodesic laminations since $\nu$ is disjoint from $C_\lambda$ and $\lambda$ is filling in $S(\lambda) \setminus C_\lambda$.
Since $A_i$ can also be regarded as an essential annulus in $H$, by Lemma \ref{esan}, each component of $\partial A_i$ is disc-dodging. Thus we have proved that  $\lambda$ is contained in   a Hausdorff limit of a sequence of disc-dodging simple closed curves on $\partial H$ as its unique minimal component.
\end{proof}

Now we are ready to prove the necessity part of Theorem \ref{main with parabolics}.

\begin{proof}[Proof of Theorem \ref{main with parabolics} (necessity)]
Consider a representation $\rho$ that has a disc-dodging parabolic or ending lamination and let us show that $\rho$ is not primitive.
By Lemma \ref{parabolic locus}, we can  assume that every parabolic curve is disc-busting. Suppose that there is a disc-dodging ending lamination $\lambda$. Since all component of $C_\lambda$ are parabolics, they are disc-busting. 
By Lemma \ref{ending lamination}, there is a sequence of disc-dodging simple closed curves $c_i$ whose Hausdorff limit contains $\lambda$ as its unique minimal component. 

Since $c_i$ is disc-dodging, there is a  meridian $m_i$ disjoint 
 from $c_i$.  As we saw in the proof of Lemma \ref{parabolic locus}, there is a free decomposition $F=A_i*B_i$ such that $c_i$ is conjugate into $A_i$.
We take a primitive closed curve $d_i$ conjugate into $B_i$.
We fix some arcs connecting a basepoint to $c_i$ and $d_i$, and regard them as elements in $\pi_1(H)=F$.
We can then consider an element of $F=\pi_1(H)$ represented as $d_i c_i^{n_i}$ which is 
 also primitive. 
 By choosing a sufficiently large $n_i$ for each $i$, we can make the closed geodesic $e_i^*$ in $\H^3/\rho(F)$
 representing $d_i c_i^{n_i}$  pass a very thin neighbourhood of the 
 closed geodesic representing $c_i$, which we denote by $c_i^*$.
 Since $c_i$ converges in the Hausdorff topology to the union of the ending lamination $\lambda$ and extra isolated leaves spiralling around $\lambda$, the closed geodesic $c_i^*$ exits every compact.
 This shows that the closed geodesics $e_i^*$ representing primitive classes 
 do not stay in a compact set.
 This implies that  $\rho$ is not primitive stable by Lemma \ref{4.1}.
\end{proof}

We note that Theorem \ref{main with parabolics} also holds when there are no parabolics, and hence Theorem \ref{main without parabolics} could be viewed as a special case of Theorem \ref{main with parabolics}. 
Indeed the necessary conditions given in \ref{main with parabolics} are automatically satisfied when there are no parabolics.
The ending lamination is minimal filling and in Masur domain in this case and such a lamination is disc-busting \cite{Ot}.
\medskip

\section{twisted $I$-bundle case}
Let us consider a special case when the non-cuspidal part $M_0$ is a twisted $I$-bundle.
We note that if $M_0$ is a twisted $I$-bundle, it has only one end.
If it is geometrically finite, then Minksy's result showed that it is primitive stable.
We are interested in the case when the end is geometrically infinite. We shall see our main theorem  implies the primitive stability for twisted $I$-bundles over non-orientable surfaces with non-empty boundaries.
We first note the following fact about ending laminations for twisted $I$-bundles.

 \begin{lem}
\label{twisted}
Let $W$ be a twisted $I$-bundle over a non-orientable surface $B$, and denote the associated $\partial I$-bundle contained in $\partial W$ by $S$.
A minimal filling lamination $\mu$ on $S$ can be an ending lamination of a hyperbolic $3$-manifold  whose non-cuspidal parts has a relative compact core  which  is homeomorphic to $(W, \partial W \setminus \Int S)$ as pared manifolds   if and only if $\mu$ is not isotopic to a double-cover of a lamination on $B$.
\end{lem}
\begin{proof}
The \lq\lq if" part was proved in Ohshika \cite{OhL}.
The \lq\lq only if" part, which is really relevant to our argument, is also well known.
We here give a brief proof.

Suppose that $\mu$ double covers a lamination $\bar \mu$ on $B$, and that $M$ is a hyperbolic 3-manifold as in the statement.
By assumption, $M_0$ has a relative compact core $C$ homeomorphic to $W$ as pared manifolds.
Then $\bar \mu$ is a Hausdorff limit of a sequence of simple closed curves $c_i$ on $B$.
We consider a pleated surface $f_i: B \rightarrow M$ realising $c_i$.
Let $B_0$ be a section of the base surface $B$ inside $C$ with respect to the $I$-bundle structure.
Since  the bundle is twisted, any surface in $M$ homotopic to $B_0$ in $M$ must intersect $B_0$ and hence the image of $f_i$ also intersects $B_0$ for every $i$.
Therefore $f_i$ cannot tend to an end because the set of such pleated surfaces is precompact (see \cite[I.5.2.18]{CEG}).
This shows that $\bar \mu$, hence also $\mu$ is realisable.
Thus we see that $\mu$ cannot be an ending lamination.
\end{proof}

In addition to this Lemma \ref{twisted}, we have the following.

\begin{lem}
\label{double covers}
For a twisted $I$-bundle $W$ as in Lemma \ref{twisted}, every filling, disc-dodging lamination in $S$ is isotopic to a double cover of a lamination on $B$.
\end{lem}
\begin{proof}
Suppose that $\lambda$ is a filling, disc-dodging lamination in $S$ and fix a transverse measure on $\lambda$.
Let $D_i$ be a compressing of $W$ with $i(\partial D_i, \lambda) \rightarrow 0$.
Let $\mu$ be a Hausdorff limit of $\{\partial D_i\}$ after passing to a subsequence.
Then $\mu$ cannot intersect $\lambda$ transversely.
Therefore the Hausdorff limit of $\partial D_i$ is  a union of $\lambda$ and isolated leaves whose ends spiral around $\lambda$.

Now, we shall analyse what kind of form $D_i$ can have.
By the standard technique of 3-dimensional topology (essentially due to Haken), we see that if a compressing disc is boundary-incompressible as a surface in $(M_0, S)$, then it is vertical, \ie the union of fibres over a proper arc in $B$.
In general, a compressing disc may be boundary-compressible, but after repeating boundary-compression finitely many times, it will become a union of finitely many vertical discs.
Therefore any compressing disc is obtained from disjoint vertical discs by band-sums.

We next turn to bound the number of bands contained in compressing discs uniformly.
If we perform a band-sum operation on two parallel vertical discs by a band contained in the $3$-ball cobounded by them, then we get a boundary-parallel disc, which contributes nothing by attaching to other discs.
The band sum of two isotopic (non-separating) discs is necessarily a 
separating disc.
Therefore, we can assume the disjoint vertical discs to start with do not contain more than two copies of a non-separating disc nor two copies of a separating disc.
This shows the number of vertical discs from which a compressing disc is constructed can be bounded by a constant depending only on the rank of $\pi_1(M_0)=F$.
Since we have to get a disc by band-sums, the graph obtained by regarding the discs as vertices and the bands as edges must be a tree.
Any tree has fewer edges than vertices.
Therefore, we also see that there is a bound on the number of bands depending only on the rank of $F$.

Recall that we have a sequence of compressing discs $D_i$ such that $\{\partial D_i\}$ converges to $\mu$ containing $\lambda$.
If we can boundary-compress each of $D_i$ after taking a subsequence, we let $D_i^1$ be one of the obtained compressing disc.
Since $i(\partial D_i, \partial D_i^1)=0$, we see that the Hausdorff limit $\mu^1$ of $\partial D_i^1$ does not intersect $\mu$ transversely, hence does not intersect $\lambda$ transversely either.
This means that $\mu^1$ is also a union of $\lambda$ and isolated leaves whose ends spiral around $\lambda$.
We repeat this operation until compressing discs in the sequence become boundary-incompressible.
As was remarked above, there is a bound 
independent of $i$ for the number of times which we can perform boundary-compression.
Therefore after finite steps, we get a sequence of boundary-incompressible, vertical compressing discs $D'_i$ such that $\{\partial D_i'\}$ converges to a geodesic lamination $\mu'$ which is a union of $\lambda$ and isolated leaves whose ends spiral around $\lambda$.

Since $D'_i$ is vertical, $\partial D_i' \cap S$ doubly covers a proper arc $a_i$ on $B$.
Therefore, we see that $\mu' \cap S$ doubly covers a lamination on $B$ which is a Hausdorff limit of $a_i$.
This implies that $\lambda$ also doubly covers a lamination on $B$.



\end{proof}

Now by Lemmata \ref{twisted} and \ref{double covers}, we see that if $M_0$ is a twisted $I$-bundle, its ending lamination always has to be disc-busting.
Theorem \ref{main with parabolics} then implies that $M_0$ is primitive stable as a representation of a free group.

For a twisted $I$-bundle $M$ over a non-orientable surface $B$ without boundary, there is similar work by Lee \cite{Lee}. Note that in this case, $\pi_1(M)$ is not a free group.
She defined an element $g\in \pi_1(M)$ to be primitive if it can be represented by a simple closed curve on $B$ and defined the primitive stability of a representation $\rho :\pi_1(M) \ra \pslc$ in a similar way as in the case of free groups.
She has shown that a representation $\rho$ is not primitive stable if and only if there exists a primitive $g$ such that
$\rho(g)$ is a parabolic.

\section*{Acknowledgements}
First of all, we want to thank to A. Lubotzky for sending us his article\cite{Lu}. His exposition gave us a great motivation for
this work. We also thank to M. Mj for invaluable discussions on Cannon-Thurston maps.
The third author would like to express his gratitude to KIAS for its hospitality during his stay when this work was done. 
Finally, we want to give our thanks to the referee for some corrections and helpful comments.

\space Woojin Jeon\\ School of Mathematics\\
 KIAS, Hoegiro 87, Dongdaemun-gu\\
     Seoul, 130-722, Korea\\
    \texttt{jwoojin\char`\@ kias.re.kr}\\

\space Inkang Kim\\
School of Mathematics\\
     KIAS, Hoegiro 87, Dongdaemun-gu\\
     Seoul, 130-722, Korea\\
     \texttt{inkang\char`\@ kias.re.kr}\\

\space Ken'ichi Ohshika\\
Department of Mathematics\\ Graduate School of Science \\Osaka University
     Toyanaka, Osaka 560-0043, Japan\\
     \texttt{Ohshika\char`\@ math.sci.osaka-u.ac.jp}\\

\space Cyril Lecuire\\
Laboratoire Emile Picard\\ Universite Paul Sabatier \\118 route de Narbonne\\
    31062 Toulouse Cedex4, France\\
     \texttt{Lecuire\char`\@math.ups-tlse.fr}\\

\end{document}